\documentclass{amsart}
\usepackage{amsfonts}
\usepackage{amsmath}
\usepackage{amssymb}
\usepackage{amsthm}
\newtheorem{theorem}{Theorem}

\newtheorem*{cor}{Corollary}

\begin{document}
\author{Mariusz Kwiatkowski, Mark Pankov}
\title{Opposite relation on dual polar spaces and half-spin Grassmann spaces}

\begin{abstract}
We characterize the collinearity (adjacency) relation of
half-spin Grassmann spaces in terms of
the relation to be opposite in the corresponding collinearity graphs.
Our characterization is
closely related with results given \cite{AVM}
and \cite{BlunckHavlicek}.
Also we show that this characterization does not hold for dual polar spaces.
\end{abstract}

\address{Department of Mathematics and Information Technology,
University of Warmia and Mazury, {\. Z}olnierska 14A, 10-561 Olsztyn, Poland}
\email{kfiecio@o2.pl, pankov@matman.uwm.edu.pl}

\maketitle

\begin{center}
{\it Dedicated to Prof. Helmut Karzel on the occasion of his 80-th birthday}
\end{center}

\section{Introduction}
Let $V$ be an $n$-dimensional vector space over a division ring.
Denote by ${\mathcal G}_{k}(V)$ the {\it Grassmannians}
consisting of all $k$-dimensional subspaces of $V$.
Two distinct elements of ${\mathcal G}_{k}(V)$ are called {\it adjacent}
if their intersection belongs to ${\mathcal G}_{k-1}(V)$
(the latter is equivalent to the fact that the sum of these subspaces is $(k+1)$-dimensional).
The cases when $k=1,n-1$ are non-interesting, since any two distinct elements of
${\mathcal G}_{k}(V)$ are adjacent if $k=1$ or $n-1$.
From this moment we suppose that $1<k<n-1$.

The {\it Grassmann graph} is the graph whose vertex set is ${\mathcal G}_{k}(V)$
and whose edges are pairs of adjacent elements.
By well-known Chow's theorem \cite{Chow},
every automorphism of this graph is induced by a semilinear isomorphism of $V$
to itself or to the dual vector space $V^{*}$,
and the second possibility can be realized only in the case when $n=2k$.
The Grassmann graph is connected and the distance between $S,U\in {\mathcal G}_{k}(V)$
is equal to
$$k-\dim(S\cap U)=\dim(S+U)-k$$
(the distance between two vertexes of a connected graph is
defined as the smallest number $i$ such that there is
a path of length $i$ connecting the vertexes);
in particular, the diameter of the Grassmann graph is finite.
Two elements of ${\mathcal G}_{k}(V)$ are called {\it opposite}
if the distance between them is equal to the diameter.
It follows from Blunck--Havlicek's results \cite{BlunckHavlicek}
(see also \cite{HavlicekPankov})
that the adjacency relation can be characterized in terms
of the relations to be opposite:
distinct $S_{1},S_{2}\in {\mathcal G}_{k}(V)$
are adjacent if and only if there exists
$S\in {\mathcal G}_{k}(V)\setminus \{S_{1},S_{2}\}$
such that every element of ${\mathcal G}_{k}(V)$ opposite
to $S$ is opposite to $S_{1}$ or $S_{2}$.
In particular, this implies that every bijective transformation of
${\mathcal G}_{k}(V)$ preserving the relation to be opposite
in both directions is an automorphism of the Grassmann graph.

In this note,
we characterize the collinearity (adjacency) relation of
half-spin Grassmann spaces in terms of
the relation to be opposite in the corresponding collinearity graph.
Also we give an example showing that this characterization
does not hold for dual polar spaces.

Usual and polar Grassmann spaces (in particular, dual polar spaces and half-spin Grassmann spaces)
are also known as the shadow spaces of buildings of type
$\textsf{A}_{n}$, $\textsf{C}_{n}$, and $\textsf{D}_{n}$.
So our considerations can be also motivated by
Abramenko--Van Maldeghem's result \cite{AVM} concerning the adjacency and opposite relations
in the chamber sets of spherical buildings.

The present note is a part of Master thesis of the first author
under supervision of the second author.
The authors thank Hendrik Van Maldeghem for useful discussion.

\section{Definitions and results}
Recall that a {\it partial linear space} is a pair $\Pi=(P, {\mathcal L})$,
where $P$ is a set of points and ${\mathcal L}$ is a family of line such that
each line contains at least two points, every point is on a certain line,
and for any distinct points there is at most one line containing them
(points connected by a line are called {\it collinear}).
For every partial linear space $\Pi=(P, {\mathcal L})$
there is the associated {\it collinearity graph} whose point set is $P$ and
whose vertexes are pairs of collinear points.

Following \cite{BuekenhoutShult} we define a {\it polar space} of finite rank
as a partial linear space $\Pi=(P, {\mathcal L})$ satisfying the axioms:
\begin{enumerate}
\item[(1)] on every line there are at least 3 points,
\item[(2)] a point is collinear with all points of a line
or with precisely one point of a line (Buekenhout--Shult's axiom),
\item[(3)] for each $p\in P$ there is a point non-collinear with $p$
(our polar space is non-degenerate),
\item[(4)] every flag consisting of singular subspaces is finite.
\end{enumerate}
Then all maximal singular subspaces are projective spaces of same finite dimension $n$,
and the number $n+1$ is known as the {\it rank} of our polar space.
The collinearity relation will be denoted by $\perp$:
we write $p\perp q$ if $p,q\in P$ are collinear and $p\not\perp q$ overwise.
Similarly, $X\perp Y$ means that every point of $X$ is collinear with all points of $Y$.
We denote by $X^{\perp}$ the set of all points $p\in P$ satisfying $p\perp X$.

A subset $F=\{p_{1}\ldots p_{2n+2}\}$ (recall that the rank of $\Pi$ is equal to $n+1$)
is called a {\it frame} of $\Pi$
if for every $p_i \in F$ there is precisely one point $p_j \in F$, $j\neq i$
such that $p_i \not \perp p_j$.
In what follows we will use the following well-know fact:
for every two singular subspaces there is a frame
whose points span both these subspaces.

For every natural $n\ge 2$ there are precisely the following two types
of rank $n$ polar spaces:
\begin{enumerate}
\item[($\textsf{C}_{n}$)]
every $(n-2)$-dimensional singular subspace is contained in at least three
maximal singular subspaces,
\item[($\textsf{D}_{n}$)]
every $(n-2) $-dimensional singular subspace is contained in precisely two
maximal singular subspaces,
\end{enumerate}
we say that a rank $n$ polar space is of type $\textsf{C}_{n}$ or $\textsf{D}_{n}$
if the corresponding case is realized.

Let $\Pi=(P,{\mathcal L})$ be a polar space of rank $n\geq 3$.
Denote by ${\mathcal G}_{k}(\Pi)$ the Grassmanian consisting of all $k$-dimensional singular subspaces.
A subset of ${\mathcal G}_{n-1}(\Pi)$ is called a {\it line } if it consists of all maximal singular subspaces
containing certain $M\in {\mathcal G}_{n-2}(\Pi)$.
The {\it dual polar space} ${\mathfrak G}_{n-1}(\Pi)$
is the partial linear space whose points are elements of ${\mathcal G}_{n-1}(\Pi)$
and whose lines are defined above.
If our polar space is of type $\textsf{D}_{n}$
the dual polar space is trivial: every line consists of precisely two points.
We say that two elements $S,U\in {\mathcal G}_{n-1} (\Pi)$ are {\it opposite}
and write $S\, {\rm op}\, U$ if the distance between them
in the collinearity graph of ${\mathfrak G}_{n-1}(\Pi)$ is maximal;
this is equivalent to the fact that $S$ and $U$ are disjoint.

Now suppose that our polar space is of type $\textsf{D}_{n}$, $n\geq 3$.
Then the Grassmannian ${\mathcal G}_{n-1}(\Pi)$
can be presented as the sum of two disjoint subsets
$${\mathcal O}_{\delta}(\Pi),\;\;\;\delta \in\{+,-\}$$
such that
the distance
$$d(S,U)=n-1-\dim(S\cap U)$$
(in the collinearity graph of ${\mathfrak G}_{n-1}(\Pi)$)
is even if $S,U$ belongs to the same ${\mathcal O}_{\delta}(\Pi)$ and odd otherwise.
These subsets are known as the {\it half-spin Grassmannians} of $\Pi$.
A subset of ${\mathcal O}_{\delta}(\Pi)$ is called a {\it line } if it consists of all elements of
${\mathcal O}_{\delta}(\Pi)$ containing certain $M\in {\mathcal G}_{n-3}(\Pi)$.
We get a partial linear space which will denoted by ${\mathfrak O}_{\delta}(\Pi)$.
In the case when $n=3$, any two distinct elements of ${\mathcal O}_{\delta}(\Pi)$
are connected by a line (their intersection is a single point)
and ${\mathfrak O}_{\delta}(\Pi)$ is a $3$-dimensional projective space.
If $n=4$ then ${\mathfrak O}_{\delta}(\Pi)$ is a polar space of type $\textsf{D}_4$.

As above, two elements $S,U\in {\mathcal O}_{\delta}(\Pi)$
are said to be {\it opposite}, $S\, {\rm op}\, U$, if the distance between them
in the collinearity graph of ${\mathfrak O}_{\delta}(\Pi)$ is maximal.
If $n$ is even then this is equivalent to the fact that
$S$ and $U$ are disjoint.
In the case when $n$ is odd, we have $S\, {\rm op}\, U$
if and only if the intersection of $S$ and $U$ is a single point.

\begin{theorem}
If $\Pi$ is of type {\rm $\textsf{D}_{n}$}, $n\ge 4$
then the following conditions are equivalent
\begin{enumerate}
\item[(1)] $S_1,S_2\in{\mathcal O}_{\delta}(\Pi)$ are collinear
points of ${\mathfrak O}_{\delta}(\Pi)$,
\item[(2)] there exists $S\in {\mathcal O}_{\delta}(\Pi)\setminus \{S_1,S_2\}$ such that
$U\, {\rm op}\, S$ implies that $U\,{\rm op}\,S_1$ or $U\, {\rm op}\, S_2$.
\end{enumerate}
\end{theorem}

\begin{cor}
Every bijective transformation of ${\mathcal O}_{\delta}(\Pi)$ preserving the relation to be opposite
is a collineation of ${\mathfrak O}_{\delta}(\Pi)$.
\end{cor}

In Section 4 we show that Theorem 1 does not hold for dual polar spaces
associated with sesquilinear forms.

\section{Proof of Theorem 1}

In this proof we will distinguish the following two cases:
\begin{enumerate}
\item[(I)] $n$ is even,
then $S\, {\rm op}\, U$ is equivalent to the fact that $S\cap U=\emptyset$,
\item[(II)]
$n$ is odd,
then $S\, {\rm op}\, U$ if and only if $S\cap U$ is a single point.
\end{enumerate}

$(1)\Longrightarrow (2)$.
Show that every point $S\ne S_{1},S_{2}$
on the line joining $S_1$ with $S_2$ is as required
(this line consists of all elements of ${\mathcal O}_{\delta}(\Pi)$ containing $S_{1}\cap S_{2}$).
Suppose that $U\,{\rm op}\, S$, but $U$ is not opposite to both $S_{1}$ and $S_{2}$.

{\it Case} (I).
In this case,
$U$ intersects $S_{1}$ and $S_{2}$ by
subspaces whose dimensions are not less than $1$.
We take lines $L_{i}\subset U\cap S_{i}$, $i=1,2$.
These lines do not intersect $S_{1}\cap S_{2}$.
Hence $S_{i}$ is spanned by $S_{1}\cap S_{2}$ and the line $L_{i}$.
The latter means that $L_1 \not\perp L_2$ which contradict the fact that
our lines are contained in $U$.

{\it Case} (II).
According our assumption,
the dimensions of $U\cap S_{1}$ and $U\cap S_{2}$
are not less than $2$.
Let $P_{i}$ be a plane contained in $U\cap S_{i}$, $i=1,2$.
The planes $P_{1},P_{2}$ both have a non-empty intersections with $S_1\cap S_2$
(because $S_1\cap S_2$ is $(n-3)$-dimensional).
Since $U\,{\rm op}\, S$, these intersections both are $0$-dimensional.
This implies the existence of lines $L_i\subset P_i$, $i=1,2$ disjoint with $S_1\cap S_2$.
As in the previous case, $L_1 \not\perp L_2$ which is impossible.

Therefore, in the both cases we have $U\,{\rm op}\, S_i$ for at least one $i\in \{1,2\}$.

$(2)\Longrightarrow (1)$.
We prove this implication in several steps.

First we establish that
{\it for every distinct collinear points $p_i\in S_i$ $(i=1,2)$
the line $p_1p_2$ intersects $S$.}

\begin{proof}[Proof in the case {\rm (I)}]
Suppose that the line $p_1p_2$ is disjoint with $S$.
There exists a maximal singular subspace $U$
containing $p_1p_2$ and opposite to $S$
(we can take any frame of $\Pi$ whose points span $S$ and the line $p_1p_2$,
the maximal singular subspace spanned by points of the frame and
opposite to $S$ is as required).
By our hypothesis, $U$ is opposite to $S_1$ or $S_2$;
this means that $p_1$ or $p_2$ is not in $U$
which contradicts to the fact that line $p_1p_2$ is in $U$.
\end{proof}

\begin{proof}[Proof in the case {\rm (II)}]
The intersection of $S_{1}$  and $S_{2}$ is not empty.
If the line $p_{1}p_{2}$ intersects $S_{1}\cap S_{2}$ then $p_{2}\in S_{1}$
and $p_{1}\in S_{2}$; thus there are the following two possibilities:
$$p_{1}p_{2}\subset S_{1}\cap S_{2}\;\mbox{ or }\;p_{1}p_{2}\cap (S_{1}\cap S_{2})=\emptyset.$$
In each of these cases, we can choose a point $p\in S_{1}\cap S_{2}$
which is not on the line $p_{1}p_{2}$ (in the first case, the dimension of $S_{1}\cap S_{2}$
is not less than $2$).

Consider the plane $P$ spanned by $p,p_{1},p_{2}$.
Assume that $P$ intersects $S$ precisely by a certain point.
Using the existence of a frame whose points span $P$ and $S$,
we construct a maximal singular subspace $U$ opposite to $S$ and containing $P$.
Then $U$ is opposite to at last one of $S_1, S_2$ which contradicts the fact that
the lines $pp_{1}$ and $pp_{2}$ are contained in $U$.

Now  suppose that $P\cap S=\emptyset$.
We choose a point $q$ from $S\cap P^{\perp}$
(this is possible since $n$ is not less than $4$)
and extend $\overline{P\cup \{ q \}}$
to a maximal singular subspace $U$ opposite to $S$
(using a frame whose points span $\overline{P\cup \{ q \}}$ and $S$).
The dimension of each $S_{i}\cap U$ is not less than $1$
which is impossible.

Therefore, $\dim (P\cap S)\ge 1$ and $P\cap S$
contains a line; this line intersects $p_{1}p_{2}$ (since $P$ is a plane).
\end{proof}

Our next step is the equalities
$$\dim (S\cap S_{i})=n-3,\;\;\;\;\;i=1,2.$$

\begin{proof}
Let us take a point $p\in S_{2}\setminus S$.
Then $S_{1}\cap p^{\perp}$ is a hyperplane of $S_{1}$ or it coincides with $S_{1}$.
Consider a line $L\subset S_{1}\cap p^{\perp}$.
Let $u,v$ be distinct points on this line.
The lines $up$ and $vp$ intersect $S$ by points $u'$ and $v'$,
respectively. Since $p\not\in S$, we have $p\ne u',v'$
and the points $u',v'$ are distinct.
The lines $L$ and $u'v'$ both are contained in the plane
$\overline{L\cup\{p\}}$, thus they have a non-empty intersection.
The inclusion $u'v'\subset S$ guarantees that $L$ intersects $S$.

So, every line $L\subset S_{1}\cap p^{\perp}$
has a non-empty intersection with $S$.
Thus $S$ intersects $S_{1}\cap p^{\perp}$ at least by a hyperplane.
The dimension of $S_{1}\cap p^{\perp}$ is not less than $n-2$
and we get
$$\dim(S\cap S_{1}\cap p^{\perp})\ge n-3$$
which implies that $S\cap S_{1}$ is $(n-3)$-dimensional.
Similarly, we show that the dimension of $S\cap S_{2}$ is equal to $n-3$.
\end{proof}

Now establish the equality $\dim S_1 \cap S_2=n-3$ which completes our proof.

\begin{proof}
Define                                      a
$$U:=(S\cap S_1)\cap (S\cap S_2).$$
Since $S\cap S_1$ and $S\cap S_2$ are $(n-3)$-dimensional subspaces of $S$,
one of the following possibilities is realized:
\begin{enumerate}
\item[(1)] $S\cap S_1=S\cap S_2$ and $U$ is $(n-3)$-dimensional,
\item[(2)] $\dim U=n-4$,
\item[(3)] $\dim U=n-5$.
\end{enumerate}
Since $U$ is contained in $S_{1}\cap S_{2}$,
the dimension of $S_{1}\cap S_{2}$ is equal to $n-3$ in the first and second cases.

Let $U$ be an $(n-5)$-dimensional subspace.
If $U$ does not coincide with $S_{1}\cap S_{2}$
then $S_{1}\cap S_{2}$ is $(n-3)$-dimensional.

Now suppose that $U=S_1\cap S_2$.
We take any line $L\subset S_1\setminus S$
and consider the singular subspace
$L^{\perp}\cap S_{2}$;
its dimension is not less than $n-3$.
Moreover, this subspace does not contain $S\cap S_{2}$.
Indeed,
$S$ is spanned by $S\cap S_{1}$ and $S\cap S_{2}$,
and the inclusion
$$S\cap S_{2}\subset L^{\perp}\cap S_{2}$$
implies that $L\perp S$; the latter is impossible, since
$S$ is a maximal singular subspace and $L\not \subset S$.

Therefore, there is a point $p\in S_{2}\setminus S$ satisfying $p\perp L$.
As above, we show that the intersection of $S$ with the plane $\overline{L\cup \{p\}}$
contains a line. This line intersects $L$ which contradicts $L\subset S_1\setminus S$.
This means that the third case can not be realized.
\end{proof}

\section{Example}
Let $V$ be a left vector space over a division ring $R$
and $\Omega: V\times V\to R$ be a non-degenerate reflexive
sesquilinear form of Witt index $n\ge 3$.
We write $\Pi=(P,{\mathcal L})$ for the associated polar space
($P$ and ${\mathcal L}$ are the sets of
$1$-dimensional and $2$-dimensional totally isotropic subspaces, respectively)
and suppose that it is of type $\textsf{C}_{n}$.
Every element of ${\mathcal G}_{n-1}(\Pi)$ can be obtained from a certain maximal singular subspace
of the form $\Omega$.

We assert that the following conditions are not equivalent
\begin{enumerate}
\item[(1)] $S_1,S_2\in{\mathcal G}_{n-1}(\Pi)$ are collinear
points of ${\mathfrak G}_{n-1}(\Pi)$,
\item[(2)] there exists $S\in {\mathcal G}_{n-1}(\Pi)\setminus \{S_1,S_2\}$ such that
$U\, {\rm op}\, S$ implies that $U\,{\rm op}\,S_1$ or $U\, {\rm op}\, S_2$.
\end{enumerate}
It is not difficult to see that (1) implies (2) (any $S\in {\mathcal G}_{n-1}(\Pi)\setminus \{S_1,S_2\}$
belonging to the line joining $S_{1}$ with $S_{2}$ is as required).
Now we show that (2) does not imply (1).

Let $p_{1},\dots,p_{n},q_{1},\dots,q_{n}$ be a frame of $\Pi$
such that $p_{i}\not\perp q_{i}$ for each $i$.
For some vectors $x_{1},\dots,x_{n},y_{1},\dots,y_{n}\in V$
we have
$$p_{i}=\langle x_{i}\rangle,\;q_{i}=\langle y_{i}\rangle\;\mbox{ and }\;
\Omega(x_{i},y_{i})=1.$$
The maximal singular subspaces of $\Pi$ associated with
the maximal totally isotropic subspaces
$$\langle x_{1},x_{2},x_{3},\dots,x_{n}\rangle,$$
$$\langle y_{1},y_{2},x_{3},\dots,x_{n}\rangle,$$
$$\langle x_{1}+y_{2},x_{2}-y_{1},x_{3},\dots,x_{n}\rangle$$
will be denoted by $S_{1},S_{2}$, and $S$ (respectively).
Their intersection is the $(n-3)$-dimensional singular subspace $N$
associated with $\langle x_{3},\dots,x_{n}\rangle$.

Now consider the line $L$ joining $\langle x_{1}+y_{2}\rangle$
with $\langle x_{2}-y_{1}\rangle$.
Every point on this line is of type
\begin{equation}
\langle (x_{1}+y_{2})+t(x_{2}-y_{1})\rangle,\;\;t\in R
\end{equation}
If $p\in p_{1}p_{2}\setminus\{p_{1}p_{2}\}$ and $q\in q_{1}q_{2}\setminus\{q_{1}q_{2}\}$ are collinear
then
$$p=\langle x_{1}+ ax_{2}\rangle\;\mbox{ and }\;q=\langle y_{2}-ay_{1} \rangle$$
for a certain scalar $a\in R$; every point on the line $p\,q$
is of type
\begin{equation}
\langle x_{1}+ ax_{2} + s(y_{2}-ay_{1})\rangle,\;\;s\in R
\end{equation}
The lines $L$ and $p\,q$ have a non-empty intersection
(because (1) coincides with (2) if $t=a$ and $s=1$).

Similarly, we establish that for any two collinear points $p\in S_{1}\setminus N$
and $q\in S_{1}\setminus N$ the line $p\,q$ intersects $S\setminus N$.
Therefore, if $U\in {\mathcal G}_{n-1}(\Pi)$ is opposite to $S$
then it is opposite to $S_{1}$ or $S_{2}$. However, $S_{1}$ and $S_{2}$ are not collinear.

\end{document}